\documentclass[a4paper]{article}

\usepackage{amsmath}
\usepackage{amssymb}
\usepackage{a4wide}
\usepackage{graphicx}
\usepackage[colorlinks=true]{hyperref}
\hypersetup{urlcolor=blue, citecolor=red}

\newtheorem{theorem}{Theorem}[section]

\newtheorem{definition}[theorem]{Definition}

\newcommand{\R}{\mathbb{R}}

\newcommand{\pt}{\partial_t}
\newcommand{\px}{\partial_x}
\newcommand{\fhi}{\varphi}

\renewcommand{\L}[1]{\mathbf{L^#1}}
\renewcommand{\u}{ {\bf u}}

\begin{document}
	
	\title{Modeling, Control and Numerics of Gas Networks}
	
	\author{
Martin Gugat\thanks{FAU Erlangen-N\"urnberg, Chair of Applied Analysis (Alexander von Humboldt-Professorship), Cauerstr. 11, 91058 Erlangen, Germany, \texttt{martin.gugat@fau.de}. 
The
	author acknowledges support by DFG TRR154 Project C03 and C05.
}
\and			
 Michael Herty\thanks{RWTH Aachen University, Institut f\"ur Geometrie
 und Praktische Mathematik, Templergraben 55,
	52074 Aachen, Germany, \texttt{herty@igpm.rwth-aachen.de}. The
	author acknowledges support by BMBF ENets 05M18PAA, 320021702/GRK2326,   333849990/IRTG-2379, 
	DFG EXC-2023  Internet of Production - 390621612, HIDSS-004, 
	DFG 18,19-1.}	
}
		
\date{\today}		
	
\maketitle

\begin{abstract}
	In this article we survey recent progress on mathematical results on gas flow in pipe networks with a special focus on 
	questions of control and stabilization. We briefly present the modeling of gas flow  and coupling conditions for flow through vertices of a network. Our main focus is on gas models for spatially one-dimensional flow governed by hyperbolic balance laws. We survey results on classical solutions as well as weak solutions.
	We present results on well--posedness,
	controllability, feedback stabilization,
	the inclusion of uncertainty in the models and numerical methods.
\end{abstract}



%%%%%%%%%%%%%%%%%%%%%%%%%%%%
\section{Introduction}

Gas transport on networks has been an active research topic for the past decade due to 
an increased demand in the sustainable use of natural resources and deregulation of energy markets. A particular focus has been on improved mathematical models of the physical flow for efficient simulation and for control purposes. 
In this article we contribute to the discussion by reviewing mathematical results on modeling, control and numerical methods for partial differential equations. Continuum mechanical models for gas flow in a single pipe are typically of the following type 
\begin{equation}
\label{eq:0} 
\partial_t \begin{pmatrix}\rho \\ \rho v \end{pmatrix} + 
\partial_x \begin{pmatrix}\rho v \\ \rho v ^2 + p(\rho) 
\end{pmatrix} = - \begin{pmatrix} 0 \\ f \; \rho v |v| + g \sin(\alpha) \rho \end{pmatrix}
\end{equation}
where $\rho$ is the gas density, $v$ the gas velocity,
$f$ is a friction parameter,
$p(\rho)$ denotes the pressure as a function of $\rho$,
$g$ is the gravitational constant and
$\alpha$ is the slope of the pipe
and where a spatially one-dimensional model is chosen due to the particularities of pipe flow with low Mach number. Those  also imply that temperature can be neglected. However, the model is capable of capturing transient phenomena driven by the need of
simulating transient gas flow. Those flow patterns appear e.g. in the case of starting up gas power plants. If one is interested in average states of the gas network simplified models might be sufficient. The model \eqref{eq:0} is posed on a single pipe and it is coupled by suitable transmission conditions to a flow model on networked pipes. The coupling typically leads to boundary conditions for the hyperbolic partial differential equation. Valves, compressor or generator stations are actuators of the control system. Those act  pointwise (not distributed) in space and they are modelled by control through the boundary conditions 
and the corresponding coupling conditions in the networked system.
Regarding  detailed flow models, numerical challenges and further results, we cannot provide a complete list of reference at this point but refer to  \cite{Benner2018,Sundar2018a,Zlotnik2015b,Brouwer2011,Bales2009,Osiadacz1984}, the references therein and the forthcoming sections. In the mathematical literature the study of hyperbolic balance laws, like \eqref{eq:0}, on metric graphs or networks has been studied over  the past decade and we refer to the  articles \cite{Bressan2014b,Canic2017a,Seo2017a,Garavello2010a} and references therein for further details. 

\setcounter{tocdepth}{1} 
\tableofcontents	
	
\section{Modeling of Gas Flow } 	

We are interested in the evolution equation (\ref{eq:1}) below 
	on a network (or metric  graph), consisting of 
	a collection of one-dimensional manifolds connected at nodes. 
	A general network is represented by
	a directed graph $\mathcal G= (\mathcal{E},\mathcal{V})$ composed by a finite number
	of edges $\mathcal{E}$ connected by vertices or junctions.
	We use a directed graph to model the network but the discussion is
	not limited to this.	For the sake of simplicity,
	we restrict the attention to a 
	%very special 
	network that is 
	composed by $n$ edges labeled by $j$ and modeled by the
	real interval $I_j:=(0,L_j)$ where we could include the case 
	$L_j=\infty$.	We remark that this simplification does not imply a loss of generality due to the fact that the hyperbolic systems \eqref{eq:1} has the property that waves propagate with finite speed. 	The density and the flux on edge $j \in \mathcal{E}$ will be  denoted $\rho_j$ and $q_j=(\rho v)_j.$ Without loss of generality the vertex is situated at  $x_*=0$ for all $j \in\mathcal{E}.$ Note that in this particular network  $v_j>0$ 
implies the flow is emerging from the vertex.  For a large, general network in order to achieve this property it might be necessary to transform the orientation of the flow at a vertex of the graph. This is achieved on the adjacent edge $j$ using the transformation $x \to -x$ and $q_j(t,x) \to - q_j(t,x)$. Note that this does not change the balance law \eqref{eq:1}.

%%%%%%%%%%%%%%%%%%%%%%%%%%%%%%%%%%%%%%%%
	
\subsection{Models for Gas Flow on Edges}

For readability we skip the edge index  in this section. Different scales might be relevant when modeling gas flow in pipe networks and we refer e.g. to \cite{Bales2009,Brouwer2011,Herty2010} for a detailed analysis of the involved scales as well as a corresponding hierarchy of suitable models. Here, we focus on a fine scale model for gas flow 
as given by the isentropic Euler equations 
\begin{equation}
\label{eq:1} 
\partial_t \begin{pmatrix}\rho \\ \rho v \end{pmatrix} + 
\partial_x \begin{pmatrix}\rho v \\ \rho v ^2 + p(\rho) 
\end{pmatrix} = - \begin{pmatrix} 0 \\ f \; \rho v |v| + g \sin(\alpha) \rho \end{pmatrix}
\end{equation}
for the gas density $\rho=\rho(t,x)$ and the gas velocity $v=v(t,x).$ The gas flux is defined by 
$q(t,x)=(\rho v)(t,x)$ and $ t\geq 0, x \in I.$ The parameters $f$ and $\alpha=\alpha(x)$ are the friction coefficient and the slope of the pipe. 
The number $g$ denotes the gravitational force. It has been shown that often it is not necessary to consider temperature variations \cite{Osiadacz1984,Osiadacz2001a}. 
Recent publication focuses on realistic pressure laws extending the isentropic choice
\begin{equation}\label{pressure-isen} p(\rho) = \kappa \rho^\gamma \end{equation}
for $\kappa>0$ and $\gamma \in [1,3]$ towards
\begin{equation}\label{pressure} p(\rho) = z(p) \rho^\gamma \end{equation}
for some polynomial $\rho \to z(\rho)$ of order at most two \cite{Grundel2015,Gugat2018c}.  
For the $z$-factor, often an affine linear model of the form
\begin{equation}
\label{zfactor}
z(p) = R \Theta \, ( 1 + \alpha \, p)
\end{equation}
is used, where $\alpha \in (-0.9, 0)$.
Equation (\ref{zfactor}) is sufficiently accurate within the network operating range,
see \cite{Gugat2018c,zbMATH06943866}.
%(see de Almeida et al.
%2014).
%see \cite{Gugat2018c}.
Often ideal gas corresponding to $\gamma=1$ is considered for control and stabilization. 
Equation \eqref{eq:1} is accompanied by initial data
\begin{equation}\label{eq:1-IC}
\rho(0,x)=\rho_0(x) \mbox{ and } q(0,x)=q_0(x),
\end{equation}
and boundary data at $x=0$ obtained through coupling conditions discussed in detail in the forthcoming section. 
In the case of $f=\alpha=0$ the system \eqref{eq:1} can be obtained formally and rigorously \cite{Lions1994c,Bouchut1999,Holle2020a} through kinetic relaxation. For $f=f(t,x,\xi) \in \R^2$ with $\xi \in\R$ the BGK model for a fixed parameter $\epsilon>0$ reads	\begin{equation}\label{kinetic}
	\partial_t f(t,x,\xi) + \partial_x \xi f(t,x,\xi) = \frac1\epsilon \left( M[f](t,x,\xi) - f(t,x,\xi) \right). 
	\end{equation}
The Maxwellian $M[f]=(M_1,M_2) \in \R^2$ is given by 
\begin{equation}
    M_1(t,x,\xi)= \chi( \rho(t,x), \xi - v(t,x) ), \; M_2(t,x,\xi)=  \chi( \rho(t,x), \xi - v(t,x) ) \left( (1-\theta) v + \theta \xi \right)
\end{equation}
where $\chi$ is a rational polynomial and $(\rho,v)$ are the moments of the kinetic distributions, i.e., 
\begin{equation}\label{moments}
\rho(t,x) = \int f_1(t,x,\xi) d\xi, \; (\rho v)(t,x) = \int f_2(t,x,\xi) d\xi.
\end{equation}
In \cite{Holle2020a} it has been shown that for fixed $\epsilon>0$ there exist a solution $f=f_\epsilon$ to equation \eqref{kinetic} subject to suitable initial and boundary conditions. Furthermore, the sequence of associated densities and fluxes
given by equation \eqref{moments} converges towards a solution to equation \eqref{eq:1} also in the presence of coupling conditions \eqref{kinetic-cplg}. For further details we refer to \cite[Theorem 2.1]{Holle2020a}. 
\par 
In engineering application and connected with discrete decisions often stationary solutions to equation \eqref{eq:1} and \eqref{pressure} are considered. The existence of steady states in the general 
case has been established in \cite{Gugat2015b,Gugat2018c}. Those stationary solutions are important for stabilization of gas flow. Furthermore, they also appear in the design of optimal controls  
due to the turnpike phenomenon. This will be discussed in detail in the section below. 
\par 
In the domain of operation of the gas pipeline networks
(low Mach numbers and large pressure) the
flow model can also modelled by a 
degenerate parabolic model as studied e.g. by 
\cite{Bamberger1977a,Leugering2018a}. Formally, the model is obtained from equation \eqref{eq:1}
by neglecting the momentum term $\rho v^2$ and the time derivative $\partial_t (\rho v).$ After
small computations one obtains a degenerate parabolic equation in the pressure $p$ in the isentropic case. 
\par 
For further models we refer to the literature, e.g. \cite{Bales2009,Brouwer2011}.

%%%%%%%%%%%%%%%%%%%%%%%%%%%%%%%%%%%%%%%%%%%%%%

\subsection{Models for the Flow Through Vertices}

The modeling of gas flow in pipes as spatially one-dimensional flow has severe implications on the modeling
of the dynamics at the vertex. In particular, the modeling solely relies on the traces of the gas density $\rho_j(t,0+)$ and 
gas flux $q_j(t,0+)$  for all adjacent edge $j$ and possibly a control input $u(t).$ In general the modeling of the vertex
is hence given by a nonlinear function $\Psi: \R^+ \times (\R^2)^n \to \R^n$ 
	\begin{equation}
	\label{eq:2}
	\Psi \left(t,  \rho_1(t,0+), q_1(t,0+), \dots, \rho_n(t,0+), q_n(t,0+) \right) = u(t) \,\; a.e.\;  t \geq 0.
	\end{equation}
Thus the one-dimensional framework allows a rather simple modeling structure given by the form of $\Psi$, but 
very fast numerical integration. While a model in \emph{one} space
	 dimension well describes the dynamics within a pipe, it hardly covers
	 geometry effects at a junction, which is clearly an intimately $3$D
	 phenomenon.  As a consequence, the literature offers several different
	 choices for a coupling or nodal  condition $\Psi$ and equation ~\eqref{eq:2}, depending on the
	 specific needs of each particular situation. 	In engineering
	 literature, the nodal conditions $\Psi$ are typically accompanied by 
	 parameters whose values are empirically justified. A numerical study
	 of one-- and two--dimensional situations can be found e.g.
	 in~\cite{K.Banda2006,Herty2007aa,Herty2013b}. 
	 \par 
	 A detailed well--posedness analysis of equation \eqref{eq:2} and \eqref{eq:1} is deferred to the next section. 
	 However, it is important to notice that the condition \eqref{eq:2} implicitly describes possible boundary conditions 
	 for the hyperbolic differential equation \eqref{eq:1}. This has implications on the modeling of suitable functions $\Psi$
	 as seen below. 
	 \par 
	 A further important aspect in the modeling of flows through vertices is the control action, here denoted by 
	 a given function $u(t)$. Those functions  may model the closure of valves or the supplied power
	 for compressor stations \cite{Herty2007}. Most examples in the literature consider explicit control actions as 
	 shown in equation \eqref{eq:2} even so an implicit dependence on $u$ is possible.  
	 Additionally, the coupling condition might dependent explicitly on time $t\geq 0$ when fatigue of material is of importance.
    \par 
	 In the following we turn to typical examples for $\Psi.$ First, 
	 consider the case of subsonic data, i.e., 
	 \begin{equation}\label{eq:subsonic}
	     \lambda_{1}( \rho_j(t,0+), q_j(t,0+) ) < 0  < \lambda_2(\rho_j(t,0+), q_j(t,0+) ). 
	 \end{equation}
    Here, $\lambda_{k}(\rho,q)$ are the $k=1,2$ eigenvalues of the Jacobian of the flux function $f$. 
    In the subsonic case  we have therefore locally in time and phase space $n$ boundary conditions
    to be determined by equation \eqref{eq:2}. A similar consideration in the supersonic case is possible and has been investigated in \cite{Gugat2017c}. In general, the boundary conditions are not explicit    but rather obtained by physical modeling considerations. \par 
    The first element of $\Psi$ typically models the conservation of mass through the vertex and 
    henceforth reads
	 \begin{equation}
	 \label{eq:5}
	 \Psi_1 ( t,  \rho_1, q_1, \dots, \rho_n, q_n ) = \sum_{j=1}^n  q_j. 
	 \end{equation}
     Except for a coupling in the gas--to--power setting \cite{Mueller2019} there is no control 
     active at the first component of $\Psi$. The other components of $\Psi$ may impose 
     different physically desirable properties. 
     \par 
     A condition typically used in the engineering community \cite{Osiadacz1984,Banda2006} 
     is to assume equal pressure $p$ at the vertex, i.e., 
	 \begin{equation}
	 \label{eq:6}
	 \Psi_j ( t,  \rho_1, q_1, \dots, \rho_n, q_n ) = p (\rho_j) - p (\rho_1) \qquad j=2, \ldots, n. 
	 \end{equation}
     An alternative condition proposed in \cite{Banda2006,Garavello2006a} is to assume 
    the continuity of the dynamic pressure or equality of momentum flux, i.e., 
	 \begin{equation}\label{eq:6b}
	 \Psi_j ( t,  \rho_1, q_1, \dots, \rho_n, q_n )
	 =
	  \left(\frac{{q_j}^2}{\rho_j} + p (\rho_j)\right)
	 -
	  \left(\frac{{q_1}^2}{\rho_1} + p (\rho_1)\right)
	 \qquad j=2, \ldots, n \,.
	 \end{equation}
	 In \cite{Garavello2006a} also geometric information has been included in the the previous two conditions 	 and an analytical comparison of qualitative properties has been conducted in 	 \cite{Colombo2009}. The conditions \eqref{eq:6} and \eqref{eq:6b} might lead to a production of energy at the vertex as observed in \cite{Reigstad2015,Mindt2019}. A condition that preserves the energy is to assume the equality of  stagnation enthalpy or equality of Bernoulli invariant
	  \begin{equation}\label{eq:6c}
	 \Psi_j ( t,  \rho_1, q_1, \dots, \rho_n, q_n )
	 =
	 \frac12 \left( \frac{q_j}{\rho_j} \right)^2 + p'(\rho_j)
	 -
	 \frac12 \left( \frac{q_1}{\rho_1} \right)^2 + p'(\rho_1)
	 \qquad j=2, \ldots, n \,.
	 \end{equation}
	 This condition implies the conservation of energy at the vertex. As noted in \cite{Holle2020}, the 
	 equation \eqref{eq:1} itself dissipates energy and this might be also a desirable property of the coupling condition. An implicit condition ensuring this property has been introduced in \cite[Definition 1]{Holle2020}. The proposed condition is implicit in the sense that only the 
	 resulting boundary values $(\rho_k,q_k)(t,0+)$ are given but {\em not} necessarily an explicit function $\Psi.$ Hence, the starting point are constant initial data $(\hat\rho_k,\hat q_k)$ for $k=1,\dots,n$ adjacent to the vertex. Then, there exist boundary values $(\rho_k,q_k)(t,0+)$ 
	 with the following properties.
	 \begin{itemize}
	     \item Mass is conserved at the junction: \begin{equation}
				\sum_{k=1}^d q_k(t,0+) =0, \, \; a.e. \; t \geq 0.
			\end{equation}
			\item There exists $\rho_*\ge 0$ such that for each fixed $k=1,\dots,n$ the boundary values $(\rho_k,q_k)(t,0+)$ are equal to the restriction to $x>0$ of the (unique) weak entropy solution of equation \eqref{eq:1} in the sense of Lax with initial condition
			\begin{equation}\label{eq:6d}
				(\rho,q)(0+,x)=
				\begin{cases}
					(\hat{\rho}_k,\hat{q}_k),\quad& x>0,\\
					(\rho_*,0),\quad& x<0.
				\end{cases}
			\end{equation}
	 \end{itemize}
	 The existence of $\rho_* > 0$ is proven in \cite[Lemma 1]{Holle2020},  the boundary values dependent continuously on $\rho_*$ \cite[Proposition 2]{Holle2020} and {\em no} assumption on subsonic initial data is required.  Further, the previous construction can be shown to be decrease entropy at the vertex for  a large class of symmetric entropies including the physical energy. 
	 The condition \eqref{eq:6d} can also be derived by using the kinetic formulation of the isentropic Euler equations \eqref{kinetic}. Coupling conditions for $f_k=f_k(t,x,\xi)$ where $f_k$ is the kinetic particle density on edge $k=1,\dots,n$ can be formulated using a coupling condition of similar type as \eqref{eq:2}, i.e., 
	 \begin{equation} 
	 \overline{\Psi}(t, f_1(t,0+,\cdot), \dots, f_n(t,0+,\cdot) )(\xi)= u(t) \; a.e. t \geq 0, \xi \in \R. 
	 \end{equation} 
	 It has been shown in \cite[Theorem 1]{Holle2020a} that among all coupling conditions $\overline{\Psi}$ that conserve the total mass, the condition that dissipates the most energy is given by 
	 \begin{equation}\label{kinetic-cplg}
	 \overline{\Psi}_k(t, f_1(t,0+,\cdot), \dots, f_n(t,0+,\cdot) ) = M(\rho_*(t), 0 ,\xi), k=1,\dots, n. 
	 \end{equation}
    Hence, in the formal limit $\epsilon \to 0$, we expect that at the vertex a state with zero velocity
    and some (unknown) density $\rho_*$ prevails. This is precisely the condition imposed by equation \eqref{eq:6d}.
     \par 
    In the uncontrolled case $u\equiv 0$ a qualitative and quantitative comparison of the previous conditions has been shown in \cite[Section 7]{Holle2020}. 
    %Here, $p(\rho)=5 \rho^2, f=\alpha=0$ and %$n=3.$ The initial data on  edge $k$ %adjacent to the vertex is %$\hat{\rho}_k=1$ for $k=1,2,3$,  $q_1=-1$ %and $q_2=q_3=\frac12.$ 
    Conditions \eqref{eq:6}, \eqref{eq:6b} lead to boundary values resulting in an energy decay of the order of $10
   ^{-2},$ while condition \eqref{eq:6c} conserves the energy and condition \eqref{eq:6d} leads to a decay of the order of $10^{-1}.$ In the case of condition \eqref{eq:6} no waves emerge from the vertex while for all other conditions except for \eqref{eq:6d} an (nonphysical) shock wave on $k=2,3$ appears. On the other hand,  rarefaction waves on edges $k=2,3$ are observed for condition \eqref{eq:6d}. 
\par 
	 Next, we turn to the modeling of controls.	The main control is due to compressor stations in gas networks. In terms of the formulation 
	 of equation \eqref{eq:2} they are modelled by considering $n=2$, equation \eqref{eq:5} and 
	 \begin{equation}
	 \label{eq:7}
	 \Psi_2(t, \rho_1, q_1, \rho_2, q_2) 
	 =
	 q_2
	 \left(
	 \left( \frac{p (\rho_2)}{p (\rho_1)} -1\right)^{(\gamma-1)/\gamma}
	 \right) 
	 \end{equation}
    for some $\gamma \in (1,3)$. Furthermore, $u(t)= (0, \Pi(t))$ where $\Pi(t)$ is the supplied compressor energy at time $t,$ 
    see e.g. \cite{Menon2005a,Herty2007e,Colombo2009b,Gugat2011a,Schmidt2012}.  This framework
	 naturally leads to various control problems, where the open--loop
	 control  has to be chosen to satisfy suitable
	 optimality criteria that will be discussed  
	 in the forthcoming sections. 
	 In contrast to the compressor control, control of valves can not add energy to the system and only one-way flow is possible through valves. 	 In \cite{Corli2019a} a mathematical model for this situation has been proposed and it amounts to prescribe a desired flow $q_*$. The control $u$ acts in such way that the valve keeps the flow at a constant value $q_*$ if possible; otherwise it is closed. Other types of valves where bi-directional flow is possible have also been considered e.g. in \cite{Corli2018c}. 
	 \par 
	 Finally and for sake of completeness, we recall that the above conditions have been partly
	 extended to the case of the full $3\times3$ system of Euler equations
	 in~\cite{Colombo2010d,Colombo2008g}.

\section{Well-Posedness of Mathematical Models For Fixed Control Action }

In this section we consider the case of a fixed  given control action $u=u(t)$ and recall well-posedness results for classical and weak solutions. 

\subsection{Classical Solutions}

Tatsien Li and his collaborators have
been very active in the study of
semi-global classical solutions to the mixed initial-boundary value problem
of one-dimensional quasilinear hyperbolic systems
and the investigation  of exact controllability in this framework,
see for example \cite{zbMATH01652441,   zbMATH07106528}.
As an example for a  quasilinear initial boundary value problem,
let us consider a density-velocity system that includes
(\ref{eq:0}) as a special case:
\begin{equation}
\label{densyvelo} 
\partial_t 
\begin{pmatrix} 
\rho  \\  v \end{pmatrix} 
+ 
\begin{pmatrix}
v & \rho  \\
\frac{(c(\rho ))^2}{\rho } & v
\end{pmatrix}
\partial_x \begin{pmatrix}  \rho  \\ 
v
\end{pmatrix} =
- 
\begin{pmatrix}
0 \\ 
F(\rho , v) 
\end{pmatrix}
\end{equation}
where the space interval is $[0,\, L]$, 
$\rho :[0, \infty)\times [0, \, L] 
\rightarrow (0, \infty)$,
$v: [0, \infty)\times [0, \, L] 
\rightarrow  {\mathbb R} $,
$c\in C^2(  (0, \infty) ; (0, \infty))$,
$F \in C^1((0, \infty) \times {\mathbb R}; {\mathbb R})$.
The initial conditions are of the form
\begin{equation}
\label{init}
\rho (0, x) = \rho _0(x), \, v(0, x) = v_0(x), \; x\in [0, L]
 \end{equation}
where $\rho _0\in C^1([0, L])$ and $v_0 \in C^1([0, L])$
are given initial data.
As an example for boundary conditions on a time
interval $[0, T]$,
let us consider the following {\sc Dirichlet} boundary conditions for subsonic flow,
where $|v|  < c(\rho )$:
\begin{equation}
\label{boundarycondi}
 v(t, 0) = u_0(t),\, \rho (t, L) = u_L(t)
 \end{equation}
where $u_0$, $u_L\in C^1([0, T])$ are
given control functions.
Classical solutions can only exist if the 
initial conditions (\ref{init}) and
the boundary conditions (\ref{boundarycondi})
are $C^1$-compatible in the sense that
\begin{equation}
\label{c0compati}
u_0(0)= v_0(0),\, u_L(0)= \rho _0(L)
\end{equation}
and
they 
satisfy the compatibility conditions
implied by
(\ref{densyvelo}), that is 
\begin{equation}
\label{c1compati}
u_0'(0)= -\frac{(c(\rho_0(0)))^2}{\rho_0(0)}
\rho _0'(0) - v_0(0) v_0'(0) -
F(\rho_0(0), v_0(0)),
\;
u_L'(0)= - v_0(L) \, \rho_0'(L) - \rho_0(L) \, v_0'(L).
\end{equation}
A typical existence result for 
semi-global  classical solutions has
the following structure:

\begin{theorem}
\label{classical}

   Let a finite (arbitrarily large) time $T>0$ be given.
   Then there exists a number $\varepsilon(T)>0$ such that
   for all initial data and boundary data 
   for which the maximal $C^1$-norm is less than $\varepsilon(T)$
   and that satisfy  the
   corresponding $C^1$-compatibility conditions
   there exists a unique  classical solution on
   the time-interval $[0, T]$, i.e. a continuously differentiable function 
   that satisfies the initial conditions,  
   the boundary conditions and the partial differential 
    equation.
    Moreover, there exists a constant $C_0(T)>0$
    such that the $C^1$-norm of the solution
    is bounded a priori by the product of
    $C_0(T)$ and the maximal $C^1$-norm of
    the initial data and the boundary data.

\end{theorem}

  The semi-global existence results are proved by
  rewriting the problem with integral equations along
  the characteristic curves
  (see \cite{zbMATH03739988}) 
  whose slopes are given by
  the eigenvalues of the system matrix.
  In the case of (\ref{densyvelo})
  the eigenvalues are $v\pm c(\rho)$. 
  Using the notation $c^2 = \frac{\partial p}{\partial \rho}$,
  the eigenvalues for (\ref{eq:0}) can also be written as 
  $\lambda_{\pm} = v \pm c$.
  
  Starting from this formulation, similar as in
  the Picard iteration a map is defined such
  that every fixed point of this map solves the initial boundary value problem.
  Then the convergence of the fixed point iteration is shown
  and the assertion follows.
  The a priori bound is shown using
  Gronwall's Lemma.
 % For sufficiently small time intervals,
  %using Gronwall's Lemma it can be shown that
  %the map is a contraction. 
  
  Note that the fixed point iteration also allows to
  consider solutions that are only required to have Lipschitz regularity, see \cite{zbMATH06892920}.
  The characteristic curves are often referred to only as 
  characteristics.  The method of characteristics is a classical method
  for the solution of  one-dimensional quasilinear hyperbolic systems.
  In order to obtain a well-posed problem,
  the boundary conditions have to be chosen according
  to the signs of the slopes of the characteristics.
  In the operation of gas networks, subsonic flow occurs,
  that is the velocity of the gas is smaller than the sound
  speed in the gas. This implies that
  one of the eigenvalues is positive in each point
  and the other one is negative in each point.
  Hence  one family of characteristic
  curves travels from the left-hand side to the right-hand side
  and the other family of characteristics  
  travels from the right-hand side to the left hand-side.
  Therefore and
  due to the structure of
  the Riemann invariants, for subsonic flow on  a single pipe, at each end
  the value of one physical variable can be prescribed   by the boundary conditions.

  For the study of control problems, it is
  often useful to work in an Hilbert space.
  This is the reason why solution of $H^2$-regularity 
  are of interest
  that are more regular than classical solutions. 
  Solutions of this type have
  been studied in \cite{zbMATH06567166}
  where an  existence result is given in  Appendix B.
  It has the same structure as Theorem \ref{classical}
  but with the $C^1$-norm replaced with the $H^2$-norm.
  In particular, the compatibility conditions remain unchanged.

%	. Klassische Lösungen, $C1,H2$. charakteristiken. %litatisen. bastin/coron. gugat/ulbrich
%
%GUGAT. Zitieren.	
%		
%FritzJohn zitieren

  Note that the existence result for semi-global classical solutions
  Theorem  \ref{classical} can be extended to
  the case of networked systems 
  that are defined on finite graphs if
  the node conditions 
  uniquely define  the necessary boundary
  input data for each adjacent pipe, see \cite{zbMATH06892920}.

\subsection{Weak Solutions}

Even for smooth initial data $(\rho_0,q_0)$ \eqref{eq:1-IC} it is known that there may exists a time $t>0$ such that solution $(\rho(t,\cdot), q(t,\cdot)$ to \eqref{eq:1} may develop discontinuities \cite{Dafermos2010b}. Hence, the notion of weak solutions has been introduced to treat solutions of lower regularity compared with solutions in the previous section. We refer to  \cite{Dafermos2010b,Bressan2009d} for details on solutions to the Cauchy problem \eqref{eq:1},\eqref{eq:1-IC} of systems of conservation and balance laws. Over several  publications those results have extended to networked 
systems \eqref{eq:1}, \eqref{eq:1-IC} and \eqref{eq:2} and refer to \cite{Bressan2014b} 
for an overview and to for $2\times 2$ hyperbolic balance laws with fixed control action to \cite{Colombo2009c}. Next, we briefly recall the basic definition of weak solution and well-posedness following \cite{Garavello2006a,Colombo2008h,Colombo2009c}. 
\par 
In order to present the notion of weak solutions for a single vertex located at $x=0$
we define $\u_j=(\rho_j,q_j)$ as density and flux on the adjacent edge $j=1,\dots,n$. Further, we define $f(\u)=(q, \frac{q^2}\rho + p(\rho) )$ with $p(\rho)$ given by 
equation \eqref{pressure-isen}. Further, we consider a source term $g(t,x,\u)=(0, - f \rho v |v| - g \sin (\alpha) \rho ).$ Then, \eqref{eq:1}, \eqref{eq:1-IC}, \eqref{eq:2}
reads 
\begin{equation}
	\label{eq:mg:system-general}
	\left\{
	\begin{array}{l}
	\partial_t \u_1 + \partial_x f (\u_1) = g (t,x,\u_1)
	\\
	\vdots
	\\
	\partial_t \u_n + \partial_x f (\u_n) = g (t,x,\u_n)
	\end{array}
	\right.
	\end{equation}
 coupled through
	the nodal condition
	\begin{equation}
	\label{eq:mg:general-nodal-condition}
	\Psi\left(\u_1\left(t,0+\right), \cdots, \u_n\left(t,0+\right)\right) = u(t). 
	\end{equation}
	Here, for every $j \in \left\{1, \ldots, n\right\}$,
	$\u_j:[0,T) \times I \to \Omega_j$, 
	$T \in (0, +\infty]$, $I=(0,\infty),$ and $\Omega_j$ is a subset of $\R^{2}.$
	We supplement~(\ref{eq:mg:system-general})
	and~(\ref{eq:mg:general-nodal-condition}) with the
	initial condition
	\begin{equation}
	\label{eq:mg:general-initial-condition}
	\left\{
	\begin{array}{ll}
	\u_1(0,x) = \u_{1,0}(x), & x>0,
	\\
	\vdots
	\\
	\u_n(0,x) = \u_{n,0}(x), & x>0,
	\end{array}
	\right.
	\end{equation}
	where, for every $j \in \left\{1, \ldots, n\right\}$,
	$\u_{j,0}:I_j \to \Omega_j$ are given functions. For brevity we introduce the
	notation
	\begin{equation}
	\label{eq:10}
	\vec\u =
	\left[
	\begin{array}{c}
	\u_1\\ \ldots\\ \u_n
	\end{array}
	\right]
	\,,\qquad
	\vec f (\u) =
	\left[
	\begin{array}{c}
	f (\u_1)\\ \ldots \\ f (\u_n)
	\end{array}
	\right]
	\,,\qquad
	\vec g (t, x, \u) =
	\left[
	\begin{array}{c}
	g(t, x, \u_1)\\ \ldots \\ g (t, x, \u_n)
	\end{array}
	\right] \,,
	\end{equation}
	and we rewrite~(\ref{eq:mg:system-general})-(\ref{eq:mg:general-nodal-condition})-(\ref{eq:mg:general-initial-condition}) in the form
	\begin{equation}
	\label{eq:4}
	\partial_t \vec \u + \partial_x \vec f (\u) = \vec g(t,x,\vec \u), \; 
	\Psi \left(\vec \u (t, 0+)\right) = u(t), \;
	\vec \u (0,x) = \vec \u_0.
	\end{equation}
%%%%%
	\begin{definition}
		\label{def:solution-to-CP}
		Fix $\hat \u = \left(\hat \u_1, \ldots, \hat \u_n\right)\in 
		\prod_{j=1}^n \Omega_j$ and $T \in ]0, +\infty]$. Assume $u \in BV(\R^+;\R^n)$. 
		A weak solution to the Cauchy problem~(\ref{eq:4})
		on $[0, T)$ is a function
		$\vec \u \in C^0 \left([0, T); \hat \u + L^1 (\R^+; \Omega^n )\right)$
		such that the following conditions hold.
		\begin{enumerate}
			\item For all $\phi \in C^\infty_c \left(]-\infty, T[ \times \R^+; \R\right)$
			and for $j \in \{1, \ldots, n\}$
			\begin{equation}
			\label{eq:weak-equality}
			\int_0^T \int_{\R^+} \left(\u_j \pt \fhi - g(t,x,\u_j) \fhi + f(\u_j) \px \fhi\right) dx\, dt
			+ \int_{\R^+} \u_{j,o} (x) \fhi(0,x) dx = 0.
			\end{equation}
			
			\item For a.e. $t \in (0,T)$,
			$\Psi \left( \vec \u(t, 0+) \right) = u(t)$.
		\end{enumerate}
		
		The weak solution $\vec \u$ is an entropy solution if for any convex entropy-entropy flux pair  $(\eta_j,Q_j)$, for all $\phi \in C^\infty_c \left(]-\infty, T[ \times \R^+; \R\right)$ and for all $j=1,\dots,n$
		\begin{equation}
			\label{eq:weak-inequality}
			\int_0^T \int_{\R^+} \left( \eta_j(\u_j) \pt \fhi - g(t,x,\u_j) D\eta_j(\u_j) \fhi + Q(\u_j) \px \fhi\right) dx\, dt
			\geq  0.
			\end{equation}
	\end{definition}
	
Well-posedness of equation \eqref{eq:4} is obtained for initial data having small total variation and under suitable assumptions on $\vec f$ and $\vec g$ as given by \cite[Assumption F and G]{Colombo2009c}. Those assumptions are precisely as in the case of the Cauchy problem. In the following we only recall the additional assumptions on $\Psi$ specific to the case of a network. We consider solutions in the space 
\begin{equation}\label{mh:weak}
    X=\left( \vec\u_* + L^1(\R^+;\Omega), u_* + L^1(\R^+; \R^n)
    \right), \; \Omega=\prod\limits_{j=1}^n \Omega_j.
\end{equation}
where the constant states $(\vec u_*, u_*)$ are such that they fulfill the coupling condition 
\begin{equation}
    \Psi(\vec u_*) = u_*.
\end{equation}
A solution is then obtained for data 
$(\vec u_0, u)= (\vec u_*, u_*) + (L^1(\R^+;\Omega), L^1(\R^+;\R^n))$ 
having small norm
\begin{equation}\label{mh:a1}
    TV(\vec \u_0) + TV(u) + \| \Psi(\vec \u_0(0+)) - u(0+) \|_{L^1(\R^+;\R^n)} \leq \delta,
\end{equation}
for some possibly small $\delta>0$. In the following we denote by $r_2(\u)$ the right eigenvector of $Df(\u)$ corresponding to the second characteristic family. Then, the assumption on $\Psi$ is as follows: 
Assume $\Psi \in C^1(\Omega;\R^n)$ such that 
\begin{equation}\label{mh:a2}
    \mbox{det} \left[ 
    D_1 \Psi(\vec \u_*) r_2(\u_{1,*}) \; 
    D_2 \Psi(\vec \u_*) r_2(\u_{2,*})
    \dots
    D_n \Psi(\vec \u_*) r_2(\u_{n,*})
    \right] \not = 0. 
\end{equation}
Under the assumptions \eqref{mh:a1} and \eqref{mh:a2} the result \cite[Theorem 2.3]{Colombo2009c} yields existence of a weak solution in the sense of Definition \ref{def:solution-to-CP}. 
\par 
Some remarks are in order. The boundary traces of $\vec u(t,\cdot)$ are well-defined, since in fact we have  $\u_j(t,\cdot) \in BV(\R^+;\Omega_j).$ 
However, the known result is only  a perturbation result around the state $(\vec \u_*, u_*).$ Theorem~2.3 \cite{Colombo2009c} also shows continuous dependence of $\vec \u$
with respect to the initial data $(\vec \u_0,u).$ This leads to an existence result for suitable optimal control problems shown in the next section. 
It is important to mention that the stated nodal conditions in the previous section {\bf all}
fulfill condition \eqref{mh:a2} 
if $\vec\u_0(x)$ is subsonic for all $x \in \R^+.$

\section{Control and Controllability }	
\subsection{Optimal Control}

In the context of weak solutions \eqref{def:solution-to-CP} existence of optimal controls has been established in \cite{Colombo2009c}. The result essentially follows by the continuous dependence of solutions $\vec \u \in X$ on $(\vec \u_0, u)$ where the space
$X$ is given by equation \eqref{mh:weak}. 

	\begin{theorem}
		\label{prop:opt}
		Let $n \in \mathbb{N}, n\geq 2$. Assume that $\vec f$
		satisfies~\textbf{(F)} at $\vec \u_*$ and $\vec G$
		satisfies~\textbf{(G)}. Fix a map $\Psi \in C^1 (\Omega;\R^n)$
		satisfying~\eqref{mh:a2} and let $u_* = \Psi(\vec \u_*)$.
	    For a fixed $\vec \u_o \in \mathcal{U}_\delta$  assume that
		\begin{eqnarray*}
			J_o:
			\left\{
			u_{\strut\vert[0,T]} \colon 
			u \in \left(u_* + \L1(\R^+;\R^n)\right)
			\mbox{ and } (\u_o,u) \in \mathcal{D}^\delta
			\right\}
			\to\R, \mbox{ and } 
			J_1: \mathcal{D}^\delta \to \R
		\end{eqnarray*}
		are non negative and lower semicontinuous with respect to the $\L1$
		norm. Then, the cost functional
		\begin{equation}
		\label{eq:J}
		\mathcal{J}(u)
		=
		J_o(u)
		+
		\int_{0}^{T} J_1\left(\vec \u(\tau,\cdot) )
		\right)d\tau
		\end{equation}
		admits a minimum on $ \left\{ u_{\strut\vert[0,T]} \colon u \in
		\left(u_* + \L1(\R^+;\R^n)\right) \mbox{ and
		}(\vec \u_o,u) \in \mathcal{D}_{0} \right\}$.
	\end{theorem}
Here, we denote by  $\vec \u(t,\cdot)$ the weak solution in the sense of 
Definition \eqref{def:solution-to-CP} with  data $(\vec \u_0, u).$ 
The definition of the sets $\mathcal{D}_0,  \mathcal{U}_\delta$ and $\mathcal{D}_\delta$ are given in \cite{Colombo2009c}. In particular, the sets $\mathcal{U}_\delta$ and $\mathcal{D}_\delta$ involve the  assumption  \eqref{mh:a1} on sufficiently small $TV$ norm.  The assumption {\bf F} and {\bf G} on flux and source term are as for the existence
of solution to a Cauchy problem and omitted here. 
\par 
In the context of gas networks a typical cost functional $\mathcal{J}$ measures the distance to a
given desired pressure $\bar{p}$  on a certain part $I_i=[x_1,x_2]$ of  
of pipe $i$. If also oscillations in the control $u(\cdot)$ should be penalised the resulting functional 
reads
\begin{equation}
    \mathcal{J}(u)= TV(u) + \int_0^T \int_{I_i} | p( \rho_i(t,x)) - \bar{p} | dxdt. 
\end{equation}
This functional fulfills the assumptions of Theorem \ref{prop:opt}. Note that 
as a possible substitute for the tracking
term for the pressure in the objective function,
also box constraints
for the pressure 
are of interest.

\par 
Since the mapping $\vec \u_0 \to \vec \u$ is non-differentiable in any $L^p, p\geq 1,$ optimality
conditions are not straightforward. An alternative 
notion of differentiability has been introduced in \cite{Bressan1998,Bressan1995c,Ulbrich2003a}. However, the extension to boundary control is still subject to active research \cite{Pfaff2017a}. 

In \cite{instantaneous}, the  instantaneous control of mixed-integer PDE-constrained gas transport problems has been studied.
Zero-one decisions occur in a natural way in the operation of gas networks
in decisions as whether to switch on or off a compressor or whether to open or close a valve.
As a first step towards the solution
of transient optimal control problems 
with decisions of this type,
an instantaneous control approach is suggested,
where in the time-discrete problem,
in each time step the control is chosen in such a way that the integrand in the objective function
for the next time step is minimized.
Approximation of the nonlinearities by piecewise linear functions
leads to large mixed integer linear  optimization problems
where a solution close to global optimality is possible.

\subsection{Controllability}
The results about the exact
boundary controllability 
of general one-dimensional quasilinear hyperbolic systems
that have been obtained in the framework of classical solutions
(see e.g. \cite{zbMATH07106528}) can be 
applied to (\ref{eq:0}).
It is typical for hyperbolic systems,
that due to the finite speed of information flow 
exact controllability is only possible
after a sufficiently large set-up time,
see \cite{zbMATH06513574}.
In a nutshell, exact controllability
is only possible after the 
flow of information from the 
boundary input 
has reached each point 
of the space interval 
on both families of characteristics.
Exact boundary controllability
on tree-like networks has also been studied in
\cite{zbMATH05969737}.

%	. Klassischen Lösungen: exakte Steuerbarkeit
%	Theorem Gugat.

%\subsubsection{Exact Controllability of
%Nodal Profiles}

In the context of gas flow
through pipelines, 
apart from the classical exact controllability, also 
the controllability of nodal profiles is 
of importance, since in the operation of gas pipeline networks,
customer satisfaction is achieved by generating the desired pressure and flow rate profiles at the the nodes
where the customers are located. 
Stated as a controllability problem,
the problem of the system operators 
is to steer the system
in such a way that at the boundary nodes,
after a finite time 
the desired nodal profiles are reached exactly 
during a certain time interval.

Note that in this notion of  controllability of nodal profiles, 
in contrast to the classical notion of exact controllability, not the full state
is prescribed at a fixed time, but
instead the boundary trace of the state is prescribed on a certain time interval.
Of course one of the variables can
simply be prescribed by the boundary conditions,
so the task is to drive simultaneously also the other variable to the values of the desired nodal profile.
This notion of 
exact controllability of nodal profiles
has been  discussed in 
\cite{nodalpofileGugatHerty2011}
in the framework of classical solutions.
The constructions in the proofs are based upon
the exchange of the roles of time and space,
that allows that the desired nodal profiles
can play the role of virtual initial conditions.
Recently, there has been a lot of
research activity in the analysis of
the exact controllability of nodal profiles
in the framework  of classical solutions
of general networked quasilinear hyperbolic systems, see
\cite{Li2010c,Gu2011,  zbMATH07233915}.
For quasilinear wave equations
it has been studied in \cite{doi:10.1002/mma.5867}
and for the Saint Venant system in
\cite{ZHUANG201934}.
%Wang2017a}.

\subsection{Feedback Stabilization}

If a desired stationary state 
is known, the problem arises,
whether the system state in the gas network can
be stabilized towards the desired state exponentially fast by suitably chosen 
boundary feedback laws.
The corresponding analysis can be  based upon suitably chosen
Lyapunov functions. 
In particular Laypunov functions
in terms of
Riemann invariants
with exponential weights have been used successfully,
see e.g.  \cite{zbMATH05679526}
in order to show the exponential decay of the 
$L^2$-norm
for a linearized system
with linear feedback laws
in terms of the Riemann invariants.
The corresponding analysis for  the quasilinear system 
on a star-shaped network is
given in \cite{zbMATH06005564}.
\par 
Extensions to the case of time delay in the feedback control have been 
presented in \cite{Gugat2011e}.

The boundary feedback stabilization of the Saint-Venant system 
by a proportional feedback control
is studied in \cite{HAYAT201952}
using a Lyapunov function in physical coordinates.
A similar analysis for general 
hyperbolic density velocity systems 
(\ref{densyvelo}) 
that also include (\ref{eq:0}) 
is presented in \cite{zbMATH06783163}.
Here linear boundary conditions of the form
\begin{equation}
\label{bastincoronfeedback}
v(t,\, 0) = k_0 \, \rho(t, \, 0),
\;
  v(t,\, L) = k_L \, \rho(t, \, L)  
\end{equation}
with real feedback parameters $k_0$, $k_L$
are considered.
In \cite{zbMATH06783163},
intervals are defined in terms of
the desired state for which the
feedback law
(\ref{bastincoronfeedback})
leads to exponential decay
of the $L^2$-norm of the distance
between the current and the desired state.
Note that in order to extend the semi-global
solutions to global solutions,
it is useful to 
work with solutions with
$H^2$ regularity
and to  show that
the $H^2$-norm of the solutions
decays exponentially fast.
This is the reason why in
\cite{zbMATH06783163, HAYAT201952}
also 
$H^2$ Lyapunov functions are presented
and the exponential decay with respect to the
$H^2$-norm is shown.

A similar analysis
with a stabilizing Neumann
feedback of the form
\[v_x(t, 0) = v^\ast_x(0) +   k_0 v_t(t, 0) ,\,
v(t, L)  = v^\ast(L)\]
for the  quasilinear wave equation 
\begin{equation}
    \label{quasilinearwave}
v_{tt} = (c^2 - v^2) \,v_{xx} - 2 \, v \,\left( v_{tx}  + (v_x)^2 \right) - 2 \, v_t \, v_x
%-\theta |v|\, v_t - \frac{3}{2} \theta v\,|v| \, v_x.
-
%\theta
%\frac{f_g}{\delta}
2 \, f
\, |v| \left( v_t +
\frac{3}{2} \, v
%\left( \frac{1}{2} \, v + v \right)
\, v_x\right)
\end{equation}
for the velocity 
that is derived from  (\ref{eq:0}) 
is given in \cite{zbMATH06742239}.
Here $v^\ast $ is a desired stationary state for 
(\ref{quasilinearwave}) and $k_0>0$ is a feedback parameter.
The Lyapunov function that is used to show the exponential decay
has the form ${\cal{L}}(t)=   {\cal{L}}_1(t) +  {\cal{L}}_2(t)$
with
\begin{eqnarray}
{\cal{L}}_1(t)
&=&
\int\limits_0^L
k\, \Big[\left(c^2-v^2 \right)\, (v-v^\ast)_x^2+ v_t^2\Big]
-2
\exp\left(- \frac{x}{L}\right)
\Big[v\, (v-v^\ast)_x^2+ v_t \, (v-v^\ast)_x\Big]\,dx,
%\int_0^L
%k\, \Big((a^2-(\bar u+u)^2)\, u_x^2+u_t^2\Big)-2
%\exp\left(- \frac{x}{L}\right)\,
%\Big((\bar u+u)\,u_x^2+u_t \, u_x\Big)\,dx,
\\
%
%\[
{\cal{L}}_2(t)
&=&
\int\limits_0^L  k\, \Big[\left(c^2-v^2\right)\,(v-v^\ast)_{xx}^2+  v_{tx}^2\Big]
-2
\exp\left(- \frac{x}{L}\right)
\Big[v  \,(v-v^\ast)_{xx}^2+ v_{tx}\,(v-v^\ast)_{xx}\Big]\,dx,
\end{eqnarray}
A typical result then yields under suitable assumptions on $k_0$ the following estimate
\begin{equation}
    {\cal{L}}(t) \leq  {\cal{L}}(0) \exp( - \nu t ), \; t \geq 0.
\end{equation}
Here, the constant $\nu>0$ is typically not known explicitly and may depend among others 
on $L$. We refer to the section on numerical results for estimates of $\nu$ in the discrete case.

The boundary feedback stabilization by proportional-integral (PI) control 
is analyzed in  \cite{zbMATH07187698}
using a Lyapunov function in physical coordinates.
The considered boundary conditions have
the differential form
\begin{equation}
\label{bastincoronfeedbackPI}
\rho(t, 0) \, v(t,\, 0) \, = Q_0(t),
\;
\rho(t, L) \,  v(t,\, L) = \kappa_L ((1 + k_L)
\rho(t, L) - Z(t)),
\,
Z'(t) = \alpha_L( \rho^\ast(L) - \rho(t, L)).
\end{equation}
where $\rho^\ast$ denotes the desired
stationary state.

\par 
Other techniques applied to stabilize linear hyperbolic balance rely on the backstepping technique and we refer e.g. to  \cite{Anfinsen2018,Krstic2008b,Krstic2008c} for further references and details. When using the backstepping technique an important aspect is the design of observers in order to determine the feedback law. State estimation and observer design have been studied in the context of system of hyperbolic equations in \cite{Anfinsen2016b,Bin2017,Anfinsen2016,Anfinsen2016a}. 
Similar techniques as  for the stabilization of gas dynamics have been used to stabilize St. Venant flow on general networks \cite{Prieur2018,Prieur2008a,Hayat2019}. The partial differential equations \eqref{eq:1} have a similar structure, except  that the pressure is given by 
\begin{equation}
    p(\rho)= \frac{g}2 \rho^2 
\end{equation}
and that there are several models available for the friction term \cite{Bastin2019a,DeHalleux2001a}. 
The system \eqref{eq:1} can also be written in quasi-linear form using the variables $(\rho, v)$ instead
of $(\rho, q).$ For classical solutions both systems are equivalent. Stabilization in terms of the variables $(\rho,v)$ has been discussed recently e.g. in 
\cite{Hayat2020arxive}.

The boundary feedback stabilization
for the degenerate parabolic model from
\cite{Bamberger1977a}
 is studied in \cite{GugatHanteJin2020}.
 Also in this contribution the
 analysis is based upon a suitably chosen
 Lyapunov function.
 The suggested feedback law has the form
 \begin{equation}
 p(t, \, 0) = p_0,  \, q(t,\, L) = \kappa_L \, p(t,\, L)
 \end{equation}
 where $p_0>0$ is a desired pressure
 value and $\kappa_L $ is a
 feedback parameter. 
 
 Note that
 in the analysis
 of the closed loop systems
 presented in this section, 
 smallness assumptions
 for the initial state
 are used.
 In the analysis, these assumptions
 imply that no shocks are generated in
 the system.
 In order to take this into account,
 in the practical application  of the feedback laws, it is important
 to keep in mind
 that it is often not clear
 whether these smallness assumptions
 are satisfied
 for the given initial states.

%%%%%%%%%%%%%%%%%%%%%%%%%%%%%%%%%%%%%%%%%%%%%%%%%%%%%%%%%%%%%%%%%%%%%%%%%%%%%%%%%%%%%%%
	
\section{Uncertainty Quantification}
In the operation of gas networks,
the treatment of uncertainties plays a
decisive role,
since customer demands 
are uncertain.
In practice,
in a procedure with several steps
the customers first buy the
option to book 
their gas consumption 
later within
a certain defined range.
Then in a second step,
precise quantities are 
nominated on a day-ahead market.
A detailed description
of the European  entry-exit 
gas market is given in
\cite{zbMATH07061088}.

In order to take  into account
 the uncertainty,
random boundary data are included in the model.
In \cite{zbMATH06657192},  a method for the computation of the probability of feasible load constellations in a stationary gas network with uncertain demand is given.
A network with a single entry and several exits with uncertain loads is studied.
The feasible flows have to satisfy  given pressure bounds in the pipes.

The numerical method is based upon
a spherical radial decomposition
that is used both for the computation
of the probabilities
and the corresponding derivatives with respect to the control.
Gradient formulae for nonlinear probabilistic constraints with non-convex quadratic forms
are presented in \cite{zbMATH07187333}. 

In order to include the
information on 
uncertainty in
the optimization problems,
probabilistic constraints of the form
\begin{equation}
\label{probabilistic}
{\mathbb P}
\left(
g(x, \omega) \geq  0
\right)
\geq 
p
\end{equation}
are useful. 
Here 
the 
 parameter
$p\in (0, 1)$ 
is a probability threshold that 
can be chosen by the
decision maker 
a priori according to his
preferences,
$x$ denotes the decision 
variable and $\omega$ is a
random variable.
The deterministic form of
(\ref{probabilistic})
is a classical inequality
constraint $g(x) \leq 0$.

For optimization problems, the
structure of the corresponding set of feasible
controls is relevant.
In general, 
for the problems with
probabilistic constraints 
this is not a convex
set. However,
in \cite{zbMATH07241833} it has been shown that 
under weak assumptions 
the feasible set is star-shaped,
which is an important result
that implies that in the spherical radial
decomposition on each ray 
at most one interval has to be considered in
the computation of the probability.
In \cite{SchusterStrauchGugatetal.2020}
this computational approach is
compared with a more general
collocation method that is based upon
kernel density estimators.

In \cite{AdelhuetteAssmannGonzalezGrandonetal.2017},
the approach with
probabilistic constraints is generalized to
a setting that also allows to take 
the dynamic case into account.
In this paper,
for a decision variable $x$,
 a set ${\cal U}$
(this could for example
be a time interval) 
a desired probability threshold
$p \in (0, 1)$
and a random variable $\omega$
probabilistic 
constraints of the form
\begin{equation}
\label{probustconstraint}
{\mathbb P}
\left(
g(x, \omega, y) \geq  0 \; \forall y
\in {\cal U}
\right)
\geq 
p
\end{equation}
are considered and referred to
as {\em probust} constraints.
This is useful for example
to define a model
where the probability
that the state satisfies the  pressure bounds throughout the time-interval
is required to be at least $p$.
Note that
(\ref{probustconstraint})
is  a stronger requirement
than the constraint
\begin{equation}
{\mathbb P}
\left(
g(x, \omega, y) \geq  0
\right)
\geq 
p\;
\forall y \in {\cal U}
\end{equation}
which does not guarantee
that for a feasible decision $x$
the constraint is satisfied
with probability $p$
uniformly for all $y\in {\cal U}$.

In \cite{zbMATH07237873},
as a step towards the treatment
of the full transient gas pipeline network flow
using probabilistic contraints,
the optimal Neumann boundary control of a vibrating string with uncertain initial data and probabilistic terminal constraints is analyzed
and a numerical method is provided.

In \cite{zbMATH06861771},
the quasilinear wave equation
(\ref{quasilinearwave})
is considered with the feedback law at
$x=0$ and uncertain boundary data
$v^\omega(t)$
at $x=L$, that is with boundary
conditions of the form 
\[v_x(t, 0) = v^\ast_x(0) +   k_0 v_t(t, 0) ,\,
v(t, L)  = v^\ast(L) + v^\omega(t).\]
It is shown that if the noise $v^\omega$
decays exponentially fast,
the $H^1$-norm of the difference
between $v$ and the desired stationary
state $v^\ast$ decays exponentially fast.

In contrast to the probabilistic
approach to robust control that
we have presented above,
the more conservative classical  approach in robust optimization
is to consider a certain range (called the uncertainty-set)
for the uncertain parameters and
optimize subject to the constraint
that all elements of the uncertainty set are feasible
(i.e. in particular the worst case scenario).
For gas networks, this approach has been
studied in \cite{doi:10.1002/net.21871}.
In order to make the approach less costly, 
in \cite{doi:10.1002/net.21871} 
a two-stage
approach is proposed.
For the two-stage model, the problem variables are classified as here-and-now variables that
have to be decided at once before the uncertainty is realized 
and wait-and-see or adjustable variables whose values can be chosen later after the uncertainty is realized.

Another aspect of uncertainty is that
also physical parameters that appear in
the pde are uncertain. 
As a contribution to this topic,
in \cite{doi:10.1080/10556788.2019.1692206}
the  problem to identify uncertain friction
parameters that vary along the pipes
is studied. 
The consequences of uncertain friction for the transport of natural gas through passive networks of pipelines
have also been studied in \cite{HeitschStrogies2019}.

%%%%%%%%%%%%%%%%%%%%%%%%%%%%%%%%%%%%%%%%%%%%%%%%%%%%%%
	
\section{Numerical Methods For Simulation and Control }	
%	. High Order Forward Discretiaztion
%	. kopplungsbedingungen im FV context
%	. Stabilization of discrete systems
\renewcommand{\u}{ {\bf u} }

The type of equation \eqref{eq:1} is a hyperbolic balance law in one-spatial 
dimension. Hence, there is a vast  literature on possible
numerical schemes to discretize system \eqref{eq:1} and we refer
e.g. to \cite{LeVeque2006} for an overview. The network structure 
imposes few particularities that will be reviewed in this section
with particular focus on finite volume schemes. For simplicity we 
also use a regular spatial grid. In order to present numerical 
discretization in compact notation we introduce
$ \u=(\rho, q)$ and $g=(0, - f \rho v |v| - g \sin(\alpha)\rho).$ 

\subsection{Discretization of Coupling Conditions for Finite Volume Schemes}
Equation \eqref{eq:1} are approximated
numerically using a finite volume method with numerical cell
size $x_{i+1}-x_i = \Delta x>0$ and time step $t^{m+1}-t^m = \Delta t$,
chosen such the CFL condition \cite{Lax2010a}
$\lambda_\text{max} \Delta t \leq \Delta x$ is satisfied. Here, 
$\lambda_\text{max}$ is the maximum absolute value of the eigenvalues
of Jacobian of $f$. The following discretization is done
for each component $\u=(\u_1,\u_2)$ separately. Within a
finite volume method
the cell average $\overline{U}_{j,i}^m$ of $\u_j$ for $j=1,2$  
of cell $i$ at time
$t^m$ is given by
\begin{equation*}
 \overline{U}_{j,i}^m:= 
 \frac{1}{ \Delta x} \int_{x_{i-\frac{1}2}}^{x_{i+\frac{1}2 }} \u_j(x,t^m) dx.
\end{equation*}
The evolution of the cell average over time $\Delta t$ is
\begin{equation}
  \label{eq:mh:generalFV}
  \overline{U}_{j,i}^{m+1} = \overline{U}_{j,i}^m - \frac{\Delta
    t}{\Delta x} \left( (
    \mathcal{F}_j)_{i+\frac{1}{2}}^m-(\mathcal{F}_j)_{i-\frac{1}{2}}^m
  \right) + \overline{G}_{j,i}^m \;,
\end{equation}
where in Godunov's method \cite{Godunov1959b}
$(\mathcal{F}_j)_{i+\frac{1}{2}}^m = \mathcal{F}_j(
\overline{U}_{j,i}^m, \; \overline{U}_{j,i+1}^m )$ denotes the
numerical flux of component $u_j(t,x), j=1,2$ 
through the boundary of the cells $i$
and $i+1$. Further,  $\overline{G}_{j,i}^n$ is an approximation to
$\frac{1}{\Delta x} \int_{t^m}^{t^{m+1}}
\int_{x_{i-\frac{1}2}}^{x_{i+\frac{1}2}} g_j(t,x,u_j(x,t)) dxdt$
obtained by a suitable quadrature rule. Within Godunov's method the
exact solution to a Riemann problem posed at the cell boundary
$i+\frac{1}2$ is used to define the numerical flux
$(\mathcal{F}_j)^m_{i+\frac{1}2}.$ Many other numerical fluxes have been proposed
and we refer to the literature \cite{LeVeque2006,Toro2009a} for further details. 
Here, we focus on Godunov's method as a basic first--order method. 
Therein, an approximation to $\u_j(t,x)$ is then obtained by the piecewise constant
reconstruction
\begin{equation}
  \u_j(t,x) = \sum_i \sum_m  \overline{U}_{j,i}^m \; 
  \chi_{ [x_{i-\frac{1}2}, x_{i+\frac{1}2}]  
    \times [t^m, t^{m+1}] }(t,x).
\end{equation}
\par 
As discussed in the previous sections the  coupling condition \eqref{eq:2}
 induces boundary conditions for equation
\eqref{eq:1}. Numerically, the corresponding discrete boundary conditions 
at $x=0+$ can be obtained by the following method. Assume
 at time $t^m$ the cell averages in the first cell $i=0$ corresponding to $x=0$ 
 of the
connected edges  are 
$\overline{U}_{j,0}^m$ for $j=1,\dots,n$.   Denote by $\sigma \to
\mathcal{L}_\kappa( \u_o, \sigma )$ the $\kappa-$th Lax curve through
the state $\u_o$ for $\kappa =1,2.$ Then, we solve for $(\sigma_1^*,
\dots, \sigma^*_n)$ using Newton's method the nonlinear system
\begin{equation}\label{sigmaFV}
  \Psi\left(t^m, \mathcal{L}_2(  \overline{U}_{1,0}^m, \sigma_1) , \dots,  
    \mathcal{L}_n(  \overline{U}_{n,0}^m, \sigma_n) \right) = u(t^m),
\end{equation}
where $u(t^m)$ denotes the discreitzed control at time $t^m.$ Under suitable
assumptions on $\Psi$ a unique solution $(\sigma_1^*,
\dots, \sigma^*_n)$ to equation \eqref{sigmaFV} exists. 
Then, a boundary value $ \overline{U}_{j,0}^{m+1}$ at time $t^{m+1}$ is 
given by equation \eqref{eq:mh:generalFV} at $i=0$ 
\begin{equation}
  \overline{U}_{j,-1}^{m} :=   \mathcal{L}_2(  \overline{U}_{j,0}^m, \sigma^*_j).
\end{equation} 
 The given construction yields a first--order approximation to the coupling
condition \eqref{eq:2}. 
\par 
For finite-volumen schemes of higher-order additional values at the boundary are required. In all recent publications \cite{Banda2016,Borsche2016a,Briani2016,Naumann2018} the additional information  required for reconstruction is obtained using a Lax-Wendroff type approach.
This amounts to differentiate condition \eqref{sigmaFV} with respect to time and 
solve the additional equations for information on the slope of a reconstructed solution $t \to \u(t,0+).$ It can be shown \cite{Banda2016} 
that this construction allows to preserve the desired order. 
\par 
Recently in   \cite{Mantri2019} the interplay of the discretization order of the  numerical flux $\mathcal{F}$ and the source term $\overline{G}$ has been investigated. In the case of spatially one-dimensional flow a numerical discretization has been proposed that allows to obtain steady-states up to machine precision even for large spatial grids $\Delta x.$ This technique has been known as well--balanced schemes for the Cauchy problem but could be extended to the case of network equations. 

%%%%%%%%%%%%%%%%%%%%%%%%%%%%%%%

\subsection{Discretization of Stabilization Problems}

The theoretical decay rates established in the previous section will be complemented by corresponding numerical results following \cite{Banda2013}. Therein, a discrete stabilization result for a first--order spatial discretization of a quasi--linear system for $\u \in \R^d$
\begin{equation}\label{quasi-linear}
    \partial_t \u + Df(\u) \partial_x \u = 0, \; \u(x,0)=\u_0(x), \; \u(t,0)=K \u(t,1)
\end{equation}
for $Df(\u)_{i,j} = \delta_{i,j} \lambda_j(\u)>0$ and $K_{ij} = \delta_{i,j} \kappa_i>0$ has been established. A  finite--volumen discretization of the quasi-linear equation \eqref{quasi-linear} 
is given by \eqref{eq:mh:generalFV} where the flux $\mathcal{F}$ is chosen as  Upwind flux  due to the fixed direction $\lambda_j(\u)>0.$ This leads to the following discretization for $j=1,\dots,d$, $i=0,\dots,N$ and $m=0,\dots,K$
\begin{equation}\label{quasi-linear-disc}
    \overline{U}_{j,i}^{m+1} = \overline{U}_{j,i}^m - \frac{\Delta
    t}{\Delta x}  \lambda_j( \overline{U}_{1,i}^m, \dots, \overline{U}_{d,i}^m )
    \left( \overline{U}_{j,i}^m  - \overline{U}_{j,i-1 }^m
    \right), \; 
\end{equation}
and for $j=1,\dots,d$, $i=0,\dots,N$ and $m \geq 0$ the initial boundary conditions are given by 
\begin{equation}
\overline{U}_{j,i }^0 = \frac{1}{\Delta x} \int_{I_i} \u_{j,0}(x) dx \mbox{ and }  
\overline{U}_{j,-1 }^m = \kappa_j \overline{U}_{j,N }^m. 
\end{equation}
We further require that initial and boundary conditions are compatible, i.e., 
\begin{equation}
    \overline{U}_{j,-1 }^0 = \kappa_j \overline{U}_{j,N }^0, 
\end{equation}
and that in the vicinity of $\u \in B_\delta(0) \subset \R^d$ the CFL condition \eqref{cfl2} holds. For given $\delta>0$,  $\Delta t$ is chosen such that
\begin{equation}\label{cfl2}
 \frac{\Delta t}{\Delta x}   \max\limits_{j=1,\dots,d} \max\limits_{ \u \in B_\delta(0) } |\lambda_j(\u) | \leq 1. 
\end{equation}
For discrete initial data $\overline{U}_{j,i }^0 \in B_\delta(0),$ having small and bounded discrete gradients discrete exponential stability of solutions to \eqref{quasi-linear-disc} 
has been established in \cite[Theorem 2]{Banda2013}. Under suitable assumptions we obtain for discrete Lyapunov function 
\begin{equation}
    L^m = \Delta x \sum\limits_{i=0}^N \sum\limits_{j=1}^d \left( 
    \overline{U}_{j,i }^m \right)^2 \exp(- \mu_j x_i ), \; m \geq 0
\end{equation}
exponential decay in time
\begin{equation} \label{decay-discrete} 
L^m \leq \exp( - \nu t^m ) L^0, \;  m \geq 0. \end{equation}
We do not review all required assumptions but recall that the exponential decay only holds provided that 
\begin{equation}
  0 <  \kappa_j^2 \leq \frac{ D_j^{\min} }{D_j^{\max}} \mbox{ and }
  \mu_j \leq \log \left( \sqrt{ \kappa_j \frac{ D_j^{\min} }{D_j^{\max}}  }^{-2} 
  \right), 
\end{equation}
where for $j=1,\dots, d$ we define $0 < D_j^{\min} := \min\limits_{ \u \in B_\delta(0) } \frac{\Delta t}{\Delta x} \lambda_j(\u) \leq \max\limits_{ \u \in B_\delta(0) } \frac{\Delta t}{\Delta x} \lambda_j(\u)  \leq D_j^{\max}\leq 1.$ The corresponding decay rate $\nu$ is given  by
\begin{equation}
    \nu =\min\limits_{j=1,\dots, d} \left( D_j^{\min} \exp( - \mu_j \Delta x ) \mu_j \frac{\Delta x}{2\Delta t}
    \right).
\end{equation}
The discrete scheme allows for explicit decay rate $\nu$ which is also  independent of the grid since $\frac{\Delta t}{\Delta x}$ is fixed. Further, several extensions e.g. to $\lambda_j(\u)<0$ and other boundary conditions are possible, see \cite{Banda2013}. In particular, in the linear case, $\lambda_j(\u)=\lambda_j$, the constant $D_j^{\min}=D_j^{\max}$. In this case $\nu = \min\limits_{j=1,\dots, d} \left( \frac12 \lambda_j \mu_j \right) $ for $\Delta x \to 0.$
\par
The results on discrete $L^2$-stability have been extended to $H^s$-norm for any $s\geq 2$ in the case of linear flux $Df(u)=A u$ and linear source terms $g(\u)=G \u$ in the recent paper
\cite{Gerster2019}.

%%%%%%%%%%%%%%%%%%%%%%%%%%%%%
\subsection{Numerical Methods for Optimal Control Problems}

As mentioned in the theoretical results, in the case of optimal boundary control problems optimality conditions 
are not straightforward to obtain. The main reason is the lack of differentiability in $L^1$ of the control to state mapping. 
A remedy in the case of systems has been introduced in \cite{Bressan1995d,Bressan2007c} and for scalar equations in \cite{Ulbrich2003a,Ulbrich2003c} 
using a novel differential. A numerical implementation of those conditions remains challenging due to the required resolution 
of the fine structure of the solution to the system \eqref{eq:1}, \eqref{eq:1-IC} and \eqref{eq:2}.  However, several approximations
of the optimality system have been studied recently in the literature and we refer to  \cite{Gugat2006a,Giles2010c,Giles2010e,Chertock2014a,Herty2014a,Herty2018a,Herty2016a,Bardos2002c,Castro2008c} for further details. 
For classical solutions, the evaluation of derivatives in the
optimal nodal control of networked hyperbolic systems
has been studied in
\cite{zbMATH05368427}.

%%%%%%%%%%%%%%%%%%%%%%%%%%%%%

\section{Open Problems}

In this section we briefly mention open problems that might be relevant for the future development of 
modeling, control and numerics for gas networks.
\par From a modeling perspective 
it would  be interesting to extend the current models in at least two directions. It has been observed that in
the operation of gas networks the quality of the injected gas varies. This leads to the
transport of a mixture of different gases and the development of suitable models for gas mixtures 
is certainly interesting and a current research topic. On the other hand,  future simulations and control of interconnected energy systems like e.g. coupled gas and electricity networks will be of  importance. From a modeling point of view this requires suitable
model couplings, from a numerical point of view this requires to treat  multi-scale phenomena due to the different involved time scales. Preliminary results in this direction have been obtained but the full control of coupled infrastructure is still at large. This might also  possibly require the  development of novel control paradigms as well as reduced models. 
\par 
From a control point of view a recent topic 
of interest 
 for the operation of large-scale gas networks is  be the turnpike phenomenon. 
Turnpike results for the complex
system dynamics can justify a combination of steady state models with dynamic models.
A survey on large time horizon control and turnpike properties for wave equations
is presented in \cite{ZUAZUA2017199}.
Due to the turnpike phenomenon,
often the control time 
can be divided into dynamic phases
at the beginning and the end of
the time interval and
a static phase between them.
In the dynamic phase one goal 
is to develop novel efficient feedback controls. A second aspect
is the development of suitable Lyapunov functions and feedback controls for  weak solutions, like BV functions, in order to treat effects
like closing valves and shock waves.
Moreover, it
would be of interest to have
feedback laws that
can stabilize the system 
for a large set of initial states,
that is to weaken the
smallness assumptions
for the initial state.
A further aspect  in the control of complex systems 
is to include reinforcement learning of dynamics that does not 
require any mathematical model. An application of those techniques towards gas networks is still an open problem.
\par 
For numerical computations novel methods to treat large-scale networks might need to be developed in order to 
obtain practical relevance. Here, tools like model-order reduction or suitable adaptive schemes might be of importance. 
See \cite{zbMATH07218214} for
a model with a 
differential-algebraic equation. 
For the efficient computation of gradients on a network an adjoint calculus is desirable. While formally this system \cite{Herty2007e} can be  derived major obstacles appear due to the non-differentiability of the control to state mapping. Here, efficient schemes would be desirable. 
\par
Finally, results obtained for energy networks might also lead to new insights for networks appearing in different transport
processes like e.g. blood flow, vehicular traffic flow or production engineering. 

\bibliographystyle{siam}
\bibliography{file}
%\bibliography{gugat2}
\end{document}